\documentclass[11pt]{article}
\usepackage[colorlinks,
            linkcolor=blue,
            anchorcolor=green,
            citecolor=red
            ]{hyperref}
\usepackage{mathrsfs,amssymb,amsfonts}
\usepackage{cite}
\usepackage{amsmath,amssymb,latexsym,color}
\usepackage{graphicx}
\usepackage[mathscr]{eucal}
\usepackage{xypic}
\newtheorem{thm}{Theorem}[section]
\newtheorem{defin}[thm]{Definition}
\newtheorem{prop}[thm]{Proposition}
\newtheorem{lemma}[thm]{Lemma}
\newtheorem{cor}[thm]{Corollary}

\newcounter{condition}[section]
\newcounter{Algorithmnumber}

\newcounter{rownumber}

\newtheorem{constr}[thm]{Construction}
\newtheorem{notation}{Notation}
\newtheorem{example}[thm]{Example}
\newtheorem{obser}[thm]{Observation}

\newtheorem{remark}{Remark}
\newenvironment{re}{\begin{remark} \rm}{\end{remark}}
\newenvironment{eg}{\begin{example} \rm}{\end{example}}
\newenvironment{ob}{\begin{obser} \rm}{\end{obser}}
\newenvironment{cons}{\begin{constr} \rm}{\end{constr}}
\newenvironment{nota}{\begin{notation} \rm}{\end{notation}}

\newcommand{\proof}{{\it Proof.\quad}}
\newcommand{\qed}{\hfill\Box\medskip}

\usepackage{CJK}
\begin{document}
\begin{CJK*}{GBK}{song}
\renewcommand{\abovewithdelims}[2]{
\genfrac{[}{]}{0pt}{}{#1}{#2}}

\title{\bf Completeness-resolvable  graphs}

\author{ Min Feng$^{1}$\quad Xuanlong Ma$^{2}$\footnote{Corresponding author.} \quad Huiling Xu$^{1}$\\
{\footnotesize   \em  $^1$School of Science, Nanjing University of Science and Technology, Nanjing, 210094, China}\\
{\footnotesize   \em $^2$School of Science, Xi'an Shiyou University, Xi'an, 710065, China}\\
}
 \date{}
 \maketitle

\begin{abstract}
  Given a connected graph $G=(V(G), E(G))$, the length of a shortest path from a vertex $u$ to a vertex $v$ is denoted by $d(u,v)$. For a proper subset $W$ of $V(G)$, let $m(W)$ be the maximum value of $d(u,v)$ as $u$ ranging over $W$ and $v$ ranging over $V(G)\setminus W$. The proper subset $W=\{w_1,\ldots,w_{|W|}\}$ is a {\em completeness-resolving set} of $G$ if
  $$
  \Psi_W: V(G)\setminus W \longrightarrow [m(W)]^{|W|},\qquad u\longmapsto (d(w_1,u),\ldots,d(w_{|W|},u))
  $$
  is a bijection, where
  $$
  [m(W)]^{|W|}=\{(a_{(1)},\ldots,a_{(|W|)})\mid 1\leq a_{(i)}\leq m(W)\text{ for each }i=1,\ldots,|W|\}.
  $$ A graph is  {\em completeness-resolvable} if it admits a completeness-resolving set.
  In this paper, we first construct the set of all completeness-resolvable graphs by using the edge coverings of some vertices in given bipartite graphs, and then establish posets on some subsets of this set by the spanning subgraph relationship. Based on each poset, we find the maximum graph and give the lower and upper  bounds for the number of edges in a minimal graph. Furthermore, minimal graphs satisfying the lower or upper bound are characterized.


\end{abstract}

\footnote{E-mail address: fgmn\_1998@163.com (M. Feng), xuanlma@mail.bnu.edu.cn (X. Ma), \\ xuhuiling@njust.edu.cn (H. Xu).}

\section{Introduction}

For a positive integer $m$, denote by $[m]$ the set of positive integers at most $m$, and for a positive integer $k$, write
$$
[m]^k=\{(a_{(1)},\ldots,a_{(k)})\mid a_{(i)}\in[m]\text{ for all }i\in [k]\}.
$$
Throughout of this paper, a graph means a finite and simple graph with at least two vertices.
Given a graph $G$, we always use $V(G)$ and $E(G)$ to denote the vertex and edge sets of $G$, respectively. The {\em order} and {\em size} of $G$ are the cardinalities of $V(G)$ and $E(G)$, respectively.
We say that $G$ is {\em connected} if for any vertices $x,y\in V(G)$, there is a path from $x$ to $y$ in $G$.
The {\em distance} between $x$ and $y$, denoted by $d(x,y)$, is the length of a shortest path from $x$ to $y$. For a proper subset $W$ of $V(G)$, write
$$
m(W)=\max\{d(w,u)\mid w\in W, u\in V(G)\setminus W\}.
$$
The proper subset $W=\{w_1,\ldots,w_{|W|}\}$ of $V(G)$ is a {\em resolving set} of $G$ if
$$
\Psi_W: V(G)\setminus W \longrightarrow [m(W)]^{|W|},\qquad u\longmapsto (d(w_1,u),\ldots,d(w_{|W|},u))
$$
is an injection.

Resolving sets of a graph were first introduced, by Slater \cite{S75} and independently, by Harary and Melter \cite{HM76} in the 1970s. Subsequently, various applications of resolving sets have appeared in the literature, as diverse as network discovery and verification \cite{Bee}, robot navigation \cite{Kh}, pharmaceutical chemistry\cite{chartrandejo}, strategies for the Mastermind game \cite{Chv},  combinatorial optimization \cite{Se}
and so on.  For an overview of resolving sets and related topics, we refer to \cite{BC11}, \cite{CH07} and \cite{jiangpoly}.

A resolving set $W$ of a connected graph $G$ is a {\em completeness-resolving set} if $\Psi_W$ is a bijection.
We say that $G$ is {\em completeness-resolvable} if $G$ admits a completeness-resolving set.
Note that every connected graph has a resolving set.

\medskip

\noindent{\bf Problem~1.} Which graphs are   completeness-resolvable?

\medskip

A vertex $x$ of a graph $G$ is {\em universal} if $x$ is adjacent to every other vertices in $G$. Denote by $\mathcal K$ the set of graphs which have a universal vertex. Let $\mathcal P$ be the set of all paths.
In this paper, we study completeness-resolvable graphs and obtain the following result.

\begin{thm}\label{main}
  Let $G$ be a connected graph. Then $G$ is completeness-resolvable if and only if $G$ is isomorphic to a graph in $\mathcal P\cup\mathcal K\cup\mathcal B\cup\mathcal C$, where $\mathcal B$ and $\mathcal C$ are as refer to Constructions~{\rm\ref{constr1}} and~{\rm\ref{constr2}}, respectively.
\end{thm}

The rest of this paper is organized as follows.

In Section~2, we first give some notions and notations, and then construct families $\mathcal B_k$ and $\mathcal C_k$ of graphs for each $k\geq 2$. The constructions depends on edge coverings of some vertices in given bipartite graphs. Actually,
the sets $\mathcal B$ and $\mathcal C$ in Theorem~\ref{main} are equal to $\bigcup_{k=2}^{+\infty}\mathcal B_k$ and $\bigcup_{k=2}^{+\infty}\mathcal C_k$, respectively.

In Section~3, Theorem~\ref{main} is proved.

In Section~4, using the spanning subgraph relationship, we establish respective posets on $\mathcal B_k$ and $\mathcal C_k$. Based on each poset, we find the maximum graph and give the lower and upper bounds for the size of a minimal graph. Furthermore, we characterize the minimal graphs satisfying the lower or upper bound.

In Section~5, we first obtain the respective ranges for diameters of graphs in $\mathcal B_k$ and $\mathcal C_k$. Then we introduce the concept of perfectness-resolvable graphs, which are closely related to completeness-resolvable graphs. Finally, we give sufficient conditions to determine a perfectness-resolvable graph, and  conclude this paper by raising a problem which graphs are perfectness-resolvable.

\section{Constructions}

We first give some notions and notations that will be used throughout of this paper.
Let $G$ be a graph. For $x\in V(G)$ and $e\in E(G)$, we say that $e$ {\em covers} $x$, or $x$ is {\em covered} by $e$, if $e$ is incident to $x$ in $G$.
For $S\subseteq V(G)$,  an {\em edge covering} of $S$, or {\em $S$-covering}, is a family $E$ of edges in $G$ such that each vertex in $S$ is covered by at least one edge in $E$.

\begin{nota}
  Given a graph $G$ and a subset $S\subseteq V(G)$, denote by $\mathcal E(G,S)$ the set of all $S$-coverings in $G$.
\end{nota}

\begin{nota}
 Let $k$ and $m$ be positive integers.

 (i) For any vector $x\in[m]^k$, denote by $x_{(i)}$ the $i$th component of $x$.

 (ii) For $I\subseteq [k]$ and $J\subseteq [m]$, write
$$
[m]^k_I(J)=\{x\in [m]^k\mid x_{(i)}\in J\text{ for all }i\in I\}.
$$
For simplify, we write $[m]^k_i(J)$, $[m]^k_I(j)$ and $[m]^k_i(j)$ instead of $[m]^k_{\{i\}}(J)$, $[m]^k_I(\{j\})$ and $[m]^k_{\{i\}}(\{j\})$, respectively.

(iii) For $S\subseteq [k]\cup[m]^k$, let $K_S$ and $\overline K_S$ denote the complete and null graphs on $S$, respectively.

(iv) Denote by $\mathcal G([k])$ and $\mathcal G([m]^k)$ the sets of all graphs with the vertex sets $[k]$ and $[m]^k$, respectively.
\end{nota}

Let $G$ be a graph. A graph $H$ is a {\em  subgraph} of $G$ if $V(H)\subseteq V(G)$ and $E(H)\subseteq E(G)$. Furthermore, the subgraph $H$ is a {\em  spanning  subgraph} of $G$
if $V(H)=V(G)$.
For $S\subseteq V(G)$, the {\em induced subgraph} of $G$ on $S$ is the graph with the vertex set $S$ such that two vertices are adjacent if and only if they are adjacent in $G$.

\begin{nota}
  Let $k$ and $m$ be positive integers. For $E\subseteq E(K_{[m]^k})$, denote by $span(E)$ the spanning subgraph of $K_{[m]^k}$ with the edge set $E$.
\end{nota}

\begin{nota}
Let $k$ and $m$ be positive integers. For $H_1\in\mathcal G([k])$ and $H_2\in\mathcal G([m]^k)$,   define $H_1\circ H_2$ as the graph on the disjoined union $[k]\sqcup[m]^k$ with the edge set $E(H_1)\cup E(H_2)\cup E(H_1,H_2)$, where
$$
E(H_1,H_2)=\{\{i,x\}\mid i\in [k],x\in [m]^k, x_{(i)}=1\}.
$$
\end{nota}

\begin{re}
For $H_1\in\mathcal G([k])$ and $H_2\in\mathcal G([m]^k)$, we have
\begin{equation}\label{eh1h2}
  |E(H_1\circ H_2)|=|E(H_1)|+|E(H_2)|+k\cdot m^{k-1}.
\end{equation}
\end{re}

Next, we construct a family $\mathcal B$ of graphs.

\begin{cons}\label{constr1}
Let $\mathcal{B}=\bigcup_{k=2}^{+\infty}\mathcal{B}_k$, where $\mathcal B_k$ is defined by the following steps.

{\rm(i)} For each $i\in [k]$, denote by $B^k_i$ the complete bipartite graph with two parts $[2]^k_i(1)$ and $[2]^k_i(2)$.

{\rm(ii)} Define
$$
\mathcal{B}_k=\{H_1\circ span(\bigcup_{i\in[k]} E_i)\mid H_1\in\mathcal{G}([k]),\; E_i\in\mathcal E(B_i^k,[2]^k_{H_1(i)}(2))\},
$$
where  $H_1(i)$ is the union of $\{i\}$ and the set of vertices adjacent to $i$ in $H_1$.
\end{cons}

Observing that $\bigcup_{i\in[k]}E(B^k_i)=E(K_{[2]^m})$, we have the following result immediately from Construction~\ref{constr1}.

\begin{prop}\label{lemmabk}
  For $k\geq 2$, pick $H_1\in\mathcal G([k])$ and $H_2\in\mathcal G([2]^k)$. Then $H_1\circ H_2\in\mathcal B_k$ if and only if the edge subset $L_i(H_2)$ is a $[2]^k_{ H_1(i)}(2)$-covering for each $i\in[k]$, where
  \begin{equation}\label{lih2}
  L_i(H_2)=E(H_2)\cap E(B^k_i).
  \end{equation}
\end{prop}

We now construct another family $\mathcal C$ of graphs.

\begin{cons}\label{constr2}
Let $\mathcal{C}=\bigcup_{k=2}^{+\infty}\mathcal{C}_k$, where $\mathcal C_k$ is constructed by the following steps.

{\rm(i)} For $i\in [k]$, let $C^k_i$ be the bipartite graph with two parts $[3]^k_i(1)$ and $[3]^k_i(2)$,
where the edge set is
$$
\{\{x,y\}\mid x\in [3]^k_i(1),y\in [3]^k_i(2), \; |x_{(t)}-y_{(t)}|\leq 1 \text{ for each } t\in[k]\setminus\{i\}\}.
$$

{\rm(ii)} For $i\in [k]$, let $D^k_i$ be the bipartite graph with two parts $[3]^k_i(2)$ and $[3]^k_i(3)$,
where the edge set is
$$
\{\{x,y\}\mid x\in [3]^k_i(2),y\in [3]^k_i(3),\;  |x_{(t)}-y_{(t)}|\leq 1 \text{ for each } t\in[k]\setminus\{i\}\}.
$$

{\rm(iii)} Define
$$\mathcal{C}_k=\{\overline K_{[k]}\circ span(\bigcup_{i\in[k]} (E_i\cup F_i))\mid E_i\in\mathcal E(C_i^k,[3]^k_i(2)),\;  F_i\in\mathcal E(D_i^k,S_i^k)\},$$
where
$S_i^k=[3]^k_{[k]}(\{2,3\})\cap [3]^k_i(3).$
\end{cons}

We get the following result immediately from the construction of $\mathcal C_k$.

\begin{prop}\label{lemmack}
  For $k\geq 2$, choose $H_2\in\mathcal G([3]^k)$. Then $\overline K_{[k]}\circ H_2\in\mathcal C_k$ if and only if the following conditions hold.

  {\rm(i)} $E(H_2)\subseteq \bigcup_{i\in[k]}(E(C^k_i)\cup E(D^k_i))$.

  {\rm(ii)} For each $i\in[k]$, the edge subset $M_i(H_2)$ is a  $[3]^k_i(2)$-covering, where
  \begin{equation}\label{mih2}
    M_i(H_2)=E(H_2)\cap E(C^k_i).
  \end{equation}

  {\rm(iii)} For each $i\in[k]$, the edge subset $N_i(H_2)$ is an $S^k_i$-covering, where
  \begin{equation}\label{nih2}
    N_i(H_2)=E(H_2)\cap E(D^k_i).
  \end{equation}
\end{prop}

\section{Proof of Theorem~\ref{main}}

In this section, we always suppose that $G$ is a connected graph.
If $W$ is a completeness-resolving set of  $G$ with $|W|=k$ and $m(W)=m$, we say that $W$ is a  {\em $(k,m)$-completeness-resolving set}, or {\em $(k,m)$-CRS} for simplify,  and $G$ is a {\em $(k,m)$-completeness-resolvable graph}, or {\em $(k,m)$-CRG} for simplify. The proof of Theorem~\ref{main} is divided in three subsections.

\subsection{$k\geq 2$ and $m=2$}

In this subsection, we determine the set of all $(k,2)$-CRGs for $k\geq2$.

\begin{lemma}\label{verifyb}
  With references to Construction~{\rm\ref{constr1}}, let $H$ be a graph in $\mathcal B_k$.

  {\rm(i)} For  $i\in[k]$ and $x\in[2]^k$, we have $d(i,x)=x_{(i)}$.

  {\rm(ii)} The graph $H$ is a $(k,2)$-CRG.
\end{lemma}
\proof (i) Write $H=H_1\circ H_2$, where
$H_1\in\mathcal G([k])$ and  $H_2\in\mathcal G([2]^k)$.
If $x_{(i)}=1$, then $d(i,x)=1=x_{(i)}$. Now suppose $x_{(i)}=2$. On one hand, since $i$ and $x$ are not adjacent in $H$, we have $d(i,x)\geq 2$. On the other hand, if $x\in [2]^k_{H_1(i)}(2)$,
by Proposition~\ref{lemmabk}, there is an edge $\{x,y\}\in L_i(H_2)$, then $y_{(i)}=1$, and so we get a path $(i,y,x)$ in $H$, which implies that $d(i,x)\leq 2$.
If $x\notin [2]^k_{H_1(i)}(2)$, then there exists a vertex $t\in H_1(i)\setminus\{i\}$ such that $x_{(t)}=1$, and so we get a path $(i,t,x)$ in $H$, which implies that $d(i,x)\leq 2$. Consequently, one has $d(i,x)=2=x_{(i)}$.

(ii) It follows from (i) that $\Phi_{[k]}(x)=x$ for each $x\in [2]^k$. Therefore, we have derived that $\Phi_{[k]}$ is a bijection, which implies that $[k]$ is a $(k,2)$-CRS, and so $H$ is a $(k,2)$-CRG, as desired.
$\qed$

\begin{prop}\label{bk}
  For $k\geq 2$, a graph $G$ is a $(k,2)$-CRG if and only if $G$ is isomorphic to a graph in $\mathcal B_k$.
\end{prop}
\proof The sufficiency holds by Lemma~\ref{verifyb} (ii). To prove the necessity, suppose that $G$ is a $(k,2)$-CRG and let $W=\{w_1,\ldots,w_k\}$ be a $(k,2)$-CRS of $G$. Then
$$
\Psi_W: V(G)\setminus W \longrightarrow [2]^{k},\qquad u\longmapsto (d(w_1,u),\ldots,d(w_{k},u))
$$
is a bijection. Define a graph $H$ on the set $[k]\cup[2]^k$ with the edge set $E_{[k]}\cup E_{[2]^k}\cup E_{[k],[2]^k}$, where
\begin{eqnarray}
E_{[k]}&=& \{\{i,j\}\mid i,j\in[k], \;  \{w_i,w_j\}\in E(G)\}, \nonumber \\
E_{[2]^k}&=& \{\{x,y\}\mid x,y\in[2]^k,\; \{\Psi_W^{-1}(x),\Psi_W^{-1}(y)\}\in E(G)\}, \nonumber  \\
E_{[k],[2]^k}&=& \{\{i,x\}\mid i\in[k],\; x\in[2]^k,\; x_{(i)}=1\}. \label{1}
\end{eqnarray}
Note that $x_{(i)}=d(w_i,\Psi_W^{-1}(x))$. It is routine to verify that
$$
\Psi: V(G) \longrightarrow [k] \cup [2]^k, \qquad  u\longmapsto
\left\{
\begin{array}{ll}
i,        &\text{if }u=w_i,\\
\Psi_W(u),&\text{if }u\in V(G)\setminus W
\end{array}\right.
$$
is an isomorphism from $G$ to $H$. Hence, graphs $G$ and $H$ are isomorphic.
Now it suffices to prove $H\in \mathcal B_k$.

Let $H_1$ and $H_2$ be the induced subgraphs of $H$ on $[k]$ and $[2]^k$, respectively. Then $H=H_1\circ H_2$ by (\ref{1}).
Let $i\in[k]$.
Pick any vertex $x\in [2]^k_{H_1(i)}(2)$. For each $t\in H_1(i)$, we have $d(w_t,\Psi_W^{-1}(x))=x_{(t)}=2$. Particularly, one gets $d(w_i,\Psi_W^{-1}(x))=2$. Hence, there is $u\in V(G)\setminus W$ such that $(w_i,u,\Psi_W^{-1}(x))$ is a path in $G$, which implies that $(i,\Psi_W(u),x)$ is a path in $H$, and so $\{\Psi_W(u),x\}\in L_i(H_2)$, where $L_i(H_2)$ is as refer to (\ref{lih2}).
From the arbitrary choice of $x$ in $[2]^k_{H_1(i)}(2)$, we have $L_i(H_2)\in\mathcal E(B_i^k, [2]^k_{H_1(i)}(2))$. It follows from Proposition~\ref{lemmabk} that $H\in \mathcal B_k$, as desired.
$\qed$

\subsection{$k\geq 2$ and $m=3$}

In this subsection, we determine the set of all $(k,3)$-CRGs for $k\geq2$.

\begin{lemma}\label{verifyc}
  With references to Construction~{\rm\ref{constr2}}, let $H$ be a graph in $\mathcal C_k$.

  {\rm(i)} For  $i\in[k]$ and $x\in[3]^k$, we have $d(i,x)=x_{(i)}$.

  {\rm(ii)} The graph $H$ is a $(k,3)$-CRG.
\end{lemma}
\proof  (i) Write $H=\overline K_{[k]}\circ H_2$, where
$H_2\in\mathcal G([3]^k)$. Note that $x_{(i)}\in\{1,2,3\}$. 

{\em Case 1.} $x_{(i)}=1$. Then $d(i,x)=1=x_{(i)}$.

{\em Case 2.} $x_{(i)}=2$.
Since $x\in [3]^k_i(2)$, from Proposition~\ref{lemmack} (ii), there is an edge $\{x,y\}\in M_i(H_2)$ for $y\in[3]^k_i(1)$,  and so we obtain a path $(i,y,x)$ in $H$, which implies that $d(i,x)=2=x_{(i)}$.

{\em Case 3.} $x_i=3$. Then $x\in [3]^k_i(3)$.
On one hand, for any $x'\in [3]^k_i(1)$, since $|x_{(i)}-x'_{(i)}|=2>1$, according to Proposition~\ref{lemmack} (i), vertices $x$ and $x'$ are not adjacent in $H_2$, which implies that $d(x',x)\geq 2$. Hence, we get $d(i,x)\geq 3$.
On the other hand, if $x\in S^k_i$, then by Proposition~\ref{lemmack} (iii), there is an edge $\{x,z\}\in N_i(H_2)$ for $z\in[3]^k_i(2)$, and further by Proposition~\ref{lemmack}~(ii), there is an edge $\{z,v\}\in M_i(H_2)$ for $v\in[3]^k_i(1)$, which implies that there is a path $(i,v,z,x)$ in $H$, and so $d(i,x)\leq 3$;
If $x\not\in S^k_i$, then there is $t\in[k]$ with $x_{(t)}=1$, which indicates that $(i,(1,\ldots,1),t,x)$ is a path in $H$, and so $d(i,x)\leq 3$.
Consequently, we get $d(i,x)=3=x_{(i)}$.

Combining all these three cases, we obtain (i).

(ii) It follows from (i) that $\Phi_{[k]}(x)=x$ for each $x\in [3]^k$. This indicates that $\Phi_{[k]}$ is a bijection,  and so $H$ is a $(k,3)$-CRG, as desired.
$\qed$

\begin{lemma}\label{null}
  Let $W$ be a $(k,3)$-CRS of a graph $G$.

   {\rm (i)} Then the induced subgraph of $G$ on $W$ is a null graph.

   {\rm (ii)} For $u,v\in V(G)\setminus W$, write $x=\Phi_W(u)$ and $y=\Phi_W(v)$.
   If $u$ and $v$ are adjacent in $G$, then $|x_{(i)}-y_{(i)}|\leq 1$ for each $i\in[k]$.
\end{lemma}
\proof Let $W=\{w_1,w_2,\ldots,w_k\}$.

(i) Suppose for the contrary that there are vertices $w_i,w_j\in W$ such that $\{w_i,w_j\}$ is an edge in $G$. Without loss of generality, assume that $\{w_i,w_j\}=\{w_1,w_2\}$. Since $\Phi_W$ is an bijection, there is a vertex $u_0\in V(G)\setminus W$ such that $\Phi_W(u_0)=(3,1,\ldots,1)$, which implies that
$$
3=d(w_1,u_0)\leq d(w_1,w_2)+d(w_2,u_0)=1+1,
$$
a contradiction.

(ii) By contradiction, suppose that there exists $i\in[k]$ such that $|x_{(i)}-y_{(i)}|\geq 2$. Noting that $x,y\in[3]^k$, we have $\{x_{(i)},y_{(i)}\}=\{1,3\}$. Without loss of generality, assume that $x_{(i)}=1$ and $y_{(i)}=3$. Then
$$
3=y_{(i)}=d(w_i,v)\leq d(w_i,u)+d(u,v)=x_{(i)}+1=2,
$$
a contradiction.
$\qed$

\begin{prop}\label{ck}
  For $k\geq 2$, a graph $G$ is a $(k,3)$-CRG if and only if $G$ is isomorphic to a graph in $\mathcal C_k$.
\end{prop}
\proof The sufficiency holds by Lemma~\ref{verifyc}~(ii).
Substituting $[3]^k$ for $[2]^k$ in the proof of Proposition~\ref{bk}, we can obtain a graph $H=H_1\circ H_2$ and an isomorphism $\Phi$ from $G$ to $H$,
where $H_1$ and $H_2$ are graphs defined on the set $[k]$ and $[3]^k$, respectively. To get the necessity, it is enough to prove $H\in \mathcal C_k$.

Noting that $\Phi^{-1}([k])$ is a $(k,3)$-CRS of $G$, we infer from Lemma~\ref{null}~(i) that $H_1=\overline K_{[k]}$.   Combining Lemma~\ref{null}~(ii) and Proposition~\ref{lemmack}, we only need to show that $M_i(H_2)\in\mathcal E(C_i^k,[3]^k_i(2))$ and $N_i(H_2)\in\mathcal E(D_i^k,S_i^k)$ for each $i\in[k]$.

Pick any vertex $x\in [3]^k_i(2)$. Then $d(i,x)=x_{(i)}=2$ by Lemma~\ref{verifyc}~(i). Since $H_1$ is null, there is a vertex $x'\in [3]^k$ such that $(i,x',x)$ is a path in $H$, which implies that $x'_{(i)}=1$, and so $\{x',x\}\in M_i(H_2)$ by Lemma~\ref{null}~(ii).
From the arbitrary choice of $x$ in $[3]^k_i(2)$, we have $M_i(H_2)\in\mathcal E(C_i^k,[3]^k_i(2))$.

Pick any vertex $y\in S_i^k$. Then $d(i,y)=y_{(i)}=3$ by Lemma~\ref{verifyc}~(i). Hence, there exist vertices $y',y''\in [k]\cup [3]^k$ such that $(i,y'',y',y)$ is a path in $H$. Since $d(t,y)=y_{(t)}\geq 2$ for each $t\in[k]$, one gets $y'\in [3]^k$.
Noting that $H_1=\overline K_{[k]}$, we have $y''\in [3]^k$ and furthermore $y''_{(i)}=1$.  It follows from Lemma~\ref{null}~(ii) that $|y''_{(i)}-y'_{(i)}|\leq 1$ and $|y'_{(i)}-y_{(i)}|\leq 1$, and so $y'_{(i)}=2$, which implies that $\{y',y\}\in N_i(H_2)$. By the arbitrary choice of $y$ in $S_i^k$, we have $N_i(H_2)\in\mathcal E(D_i^k,S_i^k)$.

The proof is completed.
$\qed$

\subsection{Proof of Theorem~\ref{main}}


\begin{prop}\label{k=1}
  A graph $G$ is a $(1,m)$-CRG if and only if $G$ is a path.
\end{prop}
\proof If $G$ is a $(1,m)$-CRG, then $G$  has a resolving set of cardinality $1$, and so $G$ is a path. Conversely, if $G$ is a path $(u_0,u_1,\ldots,u_{m})$,  then $\Phi_{\{u_0\}}(u_i)=(i)$ for $i\in[m]$, which implies that $\{u_0\}$ is a $(1,m)$-CRS, and so $G$ is a $(1,m)$-CRG, as desired.
$\qed$

\begin{prop}\label{m=1}
  A graph $G$ is a $(k,1)$-CRG if and only if the order of $G$ is $k+1$ and $G$ has a universal vertex.
\end{prop}
\proof Note that a vertex $u$ of $G$ is universal if and only if
$$
\Phi_{V(G)\setminus\{u\}}(u)=(1,\ldots,1).
$$
Hence, the desired result follows.
$\qed$

\begin{lemma}\label{lemma}
  If $G$ is a $(k,m)$-CRG with $k\geq 2$ and $m\geq 2$, then $m=2$ or $3$.
\end{lemma}
\proof By contradiction, suppose $m\geq 4$. Let $W=\{w_1,w_2,\ldots,w_k\}$ be a $(k,m)$-CRS of $G$.
Since $\Phi_W$ is a bijection, there exist vertices $u$ and $v$ in $V(G)\setminus W$ such that $\Phi_W(u)=(1,1,\ldots,1)$ and $\Phi_W(v)=(4,1,\ldots,1)$, which implies that
$$
4=d(w_1,v)\leq d(w_1,u)+d(u,w_2)+d(w_2,v)=1+1+1,
$$
a contradiction.
$\qed$

 Theorem~\ref{main} follows from Lemma~\ref{lemma} and Propositions~\ref{bk},~\ref{ck},~\ref{k=1} and~\ref{m=1}.

\section{Partially ordered sets}

In Theorem~\ref{main}, graphs in $\mathcal P\cup\mathcal K$ are well-understood, however, graphs in $\mathcal B\cup\mathcal C$ are non-intuitive. Now the letter $k$ is always used to denote a given positive integer at least $2$.
Note that graphs in $\mathcal B_k$ (resp. $ \mathcal C_k$) have the same vertex set $[2]^k$ (resp. $[3]^k$).

\begin{nota}
For graphs $G_1$ and $G_2$ with the same vertex set, define $G_1\preceq G_2$ if $G_1$ is a spanning subgraph of $G_2$. Write $G_1\prec G_2$ if $G_1\preceq G_2$ and $G_1\neq G_2$.
\end{nota}

A {\em partially ordered set}, or {\em poset}, is  an ordered pair $(A,\leq)$ such that $\leq$ is a reflexive, antisymmetric and transitive binary relation on the set $A$.  An element $x$ is {\em maximum} in $A$ if $a\leq x$ for each $a\in A$. An element $y$ is {\em minimal} in $A$ if $a\leq y$ implies $a=y$. Observe that $(\mathcal B_k,\preceq)$ and $(\mathcal C_k,\preceq)$ are posets on graphs.

In Subsection~4.1 (resp. 4.2), based on the poset $(\mathcal B_k,\preceq)$ (resp. $(\mathcal C_k,\preceq)$), we first obtain the maximum graph in $\mathcal B_k$ (resp. $\mathcal C_k$), and turn the problem of characterizing graphs in $\mathcal B_k$ (resp. $\mathcal C_k$) into characterizing minimal graphs in $\mathcal B_k$ (resp. $\mathcal C_k$). Then we investigate minimal graphs in $\mathcal B_k$ (resp. $\mathcal C_k$).

In Subsection 4.1, noting from Construction~\ref{constr1} that the first factor of a graph in $\mathcal B_k$ is an arbitrary graph in $\mathcal G([k])$, we fix $H_1\in\mathcal G([k])$ and give a definition of an {\em $H_1$-minimal} graph.
Then the size of an  $H_1$-minimal graph is bounded by using the degrees of vertices in $H_1$, and furthermore, we describe the $H_1$-minimal graphs satisfying the lower or upper bound. As examples, we apply the results to study $K_{[k]}$-minimal graphs and  $\overline K_{[k]}$-minimal graphs.

In Subsection~4.2, noting from Construction~\ref{constr2} that the first factor of a graph in $\mathcal C_k$ is the null graph in $\mathcal G([k])$, we say that $H_2$ is {\em $k$-minimal} if $\overline K_{[k]}\circ H_2$ is a minimal graph in $\mathcal C_k$.
We establish the lower and upper bounds for the size of a $k$-minimal graph, and characterize all $k$-minimal graphs satisfying the lower or upper bound.

\subsection{A poset on $\mathcal B_k$}

In this subsection, we focus on the poset $(\mathcal B_k,\preceq)$, and always suppose that $H_1\in\mathcal G([k])$ and $H_2\in\mathcal G([2]^k)$.

\begin{thm}\label{b adding edges}
  Suppose that $H_1\circ H_2\in \mathcal B_k$. If $H_1\preceq H_1'\preceq K_{[k]}$ and $H_2\preceq H_2'\preceq K_{[2]^k}$, then $H_1'\circ H_2'\in\mathcal B_k$.
\end{thm}
\proof For each $i\in[k]$, we get
$[2]^k_{H_1'(i)}(2)\subseteq [2]^k_{H_1(i)}(2)$ and $
L_i(H_2)\subseteq L_i(H_2')$, where $L_i(H_2)$ is as refer to (\ref{lih2}).
Hence, the desired result follows by Proposition~\ref{lemmabk}.
$\qed$

Let $G_1$ and $G_2$ be graphs with the same vertex set $V$. The {\em union} $G_1\cup G_2$ is the graph with the vertex set $V$ and the edge set $E(G_1)\cup E(G_2)$.

\begin{cor}\label{unionb}
 For $G_1,G_2\in\mathcal B_k$, we have $G_1\cup G_2\in \mathcal B_k$.
Particularly, we have
 $$\bigcup_{G\in \mathcal B_k}G=K_{[k]}\circ K_{[2]^k}\in\mathcal B_k.$$
\end{cor}
\proof Let $G_1=H_1\circ H_2$ and $G_2=H_1'\circ H_2'$. Then $G_1\cup G_2=(H_1\cup H_1')\circ (H_2\cup H_2')$.
Hence, the desired result follows from
Theorem~\ref{b adding edges}.
$\qed$

A {\em join-semilattice}  is a poset in which any pair of elements has the least upper bound.
Note that each finite join-semilattice has a unique maximum element. The following result is immediate from Corollary~\ref{unionb}.

\begin{prop}
  The poset $(\mathcal B_k,\preceq)$ is a finite join-semilattice with the maximum graph $K_{[k]}\circ K_{[2]^k}$.
\end{prop}

From Theorem~\ref{b adding edges}, it is nature to put forward the following notion. We say that $H_2$ is {\em $H_1$-minimal} if $H_1\circ H_2\in\mathcal B_k$ and $H_1\circ H_2'\not\in\mathcal B_k$ for $H'_2\prec H_2$.

\begin{re}\label{remark1}
  If $H_1\circ H_2$ is a minimal graph in $\mathcal B_k$, then $H_2$ is $H_1$-minimal. However, the converse is not true.
\end{re}

We get the following result from Proposition~\ref{lemmabk}.

\begin{lemma}\label{minimal bk}
  A graph $H_2$ is $H_1$-minimal if and only if the following conditions hold.

  {\rm(i)} For each $i\in[k]$, the edge subset $L_i(H_2)$, as refer to {\rm(\ref{lih2})}, is a $[2]^k_{H_1(i)}$-covering.

  {\rm(ii)} For each $e\in E(H_2)$, there exists a vertex $i\in[k]$ such that $L_i(H_2)\setminus\{e\}$ is not a $[2]^k_{H_1(i)}$-covering.
\end{lemma}

\begin{nota}
  Let $H_1\in\mathcal G([k])$ and $H_2\in \mathcal G([2]^k)$.

(i) For $x\in[2]^k$, let
$J_x=\{i\in[k]\mid x\in[2]^k_{H_1(i)}(2)\}.$

(ii) For $x\in[2]^k$ and $e\in E(H_2)$, write
$$
I_x(e)=\{i\in[k]\mid x\in[2]^k_{H_1(i)}(2),\; e\in L_i(H_2),\; e\text{ covers }x\},
$$
and furthermore, let
$$
I_x=\bigcup_{f\in E(H_2)}I_x(f)\qquad \text{ and }\qquad \tilde I_x(e)=I_x\setminus (\bigcup_{f\in E(H_2)\setminus\{e\}} I_x(f)).
$$
\end{nota}

\begin{lemma}\label{minimal bk2}
   A graph $H_2$ is $H_1$-minimal if and only if the following conditions hold.

  {\rm(i)}  For each $x\in[2]^k$, we have $I_x=J_x$.

  {\rm(ii)}  For each $e\in E(H_2)$, there is a vertex $x\in[2]^k$ such that $\tilde I_x(e)\neq\emptyset$.
\end{lemma}
\proof Note that $I_x\subseteq J_x$. Observing that $i\in J_x$ if and only if $x\in [2]^k_{H_1(i)}(2)$, we obtain (i) if and only if the condition (i) in Lemma~\ref{minimal bk} holds.

Now suppose that the condition (i) in Lemma~\ref{minimal bk} holds.
Then $i\in \tilde I_x(e)$ if and only if any edge in $L_i(H_2)\setminus\{e\}$ does not cover $x$. Hence, we obtain (ii) if and only if the condition (ii) in Lemma~\ref{minimal bk} holds.

Consequently, the desired result follows from Lemma~\ref{minimal bk}.
$\qed$

Actually, given a graph $H_1$, the $H_1$-minimal graph $H_2$ is not unique.
We give lower and upper bounds for the size of $H_2$. The {\em degree} of a vertex in a graph  is the number of edges  covering this vertex in the graph.

\begin{thm}\label{enbounds12}
If $H_2$ is $H_1$-minimal, then
  \begin{equation*}
  2^{k-\min\{d_i\mid i\in[k]\}-1} \leq |E(H_2)| \leq \sum_{i=1}^k 2^{k-d_i-1},
  \end{equation*}
  where $d_i$ is the degree of the vertex $i$ in $H_1$.
\end{thm}
\proof For each $i\in[k]$, write $E_i=L_i(H_2)$ and denote by $E_i'$ the set of edges $e$ such that $L_i(H_2)\setminus\{e\}$ is not a $[2]^k_{H_1(i)}$-covering.
By Lemma~\ref{minimal bk}~(ii), we have
$$
  \max\{|E_i|\mid i\in[k]\}\leq |\bigcup_{i\in[k]}E_i|=|E(H_2)|=|\bigcup_{i\in[k]}E_i'| \leq \sum_{i\in[k]}|E'_i|.
$$
Noting that $|E_i'|\leq 2^{k-|H_1(i)|}=2^{k-d_i-1}\leq |E_i|$, we get the desired result.
$\qed$

\begin{re}
  If $H_1\circ H_2$ is a minimal graph in $\mathcal B_k$, then combining the equation (\ref{eh1h2}), Remark~\ref{remark1} and Theorem~\ref{enbounds12}, we get
  $$
  k\cdot 2^{k-1}+\frac{1}{2}\sum_{i=1}^kd_i+2^{k-\min\{d_i\mid i\in[k]\}-1}
  \leq |E(H_1\circ H_2)|\leq k\cdot 2^{k-1}+\frac{1}{2}\sum_{i=1}^k(2^{k-d_i}+d_i).
  $$
\end{re}

The {\em minimum degree} of a graph is the minimum value of  degrees of all vertices in this graph. We use the following result to characterize the bounds in Theorem~\ref{enbounds12}.

\begin{cor}\label{tightcondition}
  Suppose that $H_2$ is $H_1$-minimal.

  {\rm(i)}  The lower bound in Theorem~{\rm\ref{enbounds12}} is attained if and only if there exists a vertex $i$ with minimum degree in $H_1$ such that  $L_i(H_2)=E(H_2)$.

  {\rm(ii)} The upper bound in Theorem~{\rm\ref{enbounds12}} is attained if and only if for any distinct vertices $x,y\in [2]^k$ and any distinct edges $e,f\in E(H_2)$, we have {\rm(a)} $|I_x(e)|\leq 1$, {\rm(b)} either $I_x(e)=\emptyset$ or $I_y(e)=\emptyset$, {\rm(c)} $I_{x}(e)\cap I_x(f)=\emptyset$.
\end{cor}
\proof (i) It is immediate from the proof of the lower bound in Theorem~\ref{enbounds12}.

(ii) Suppose that (a), (b) and (c) hold. Observing that $\tilde I_x(e)\subseteq I_x(e)$, we infer from  Lemma~\ref{minimal bk2} (ii) that
\begin{eqnarray*}
  E(H_2)&=&\bigcup_{x\in[2]^k}\{e\in E(H_2) \mid \tilde I_x(e)\neq\emptyset\}=\bigcup_{i\in[k]}\;\bigcup_{x\in [2]^k_{H_1(i)}(2)}\{e\in E(H_2) \mid i\in I_x(e)\} \\
  &=& \bigcup_{i\in[k]}\{e\in E(H_2) \mid i\in I_x(e)\text{ for a certain vertex }x\in [2]^k_{H_1(i)}(2)\},
\end{eqnarray*}
 and so $|E(H_2)|=\sum_{i\in[k]}|[2]^k_{H_1(i)}(2)|$,  attaining the upper bound in Theorem~\ref{enbounds12}.

In the following, suppose that $|E(H_2)|$ attains the upper bound in Theorem~\ref{enbounds12}.
With references to the proof of Theorem~\ref{enbounds12},  the following conditions hold.

\medskip

(C1)   $|E_i'|=2^{k-d_i-1}$ for each $i\in[k]$.

(C2) $E_i'\cap E_j'=\emptyset$ for any distinct vertices $i$ and  $j$ in $[k]$.

\medskip

\noindent By (C1),  the following condition holds for each $i\in[k]$.

\medskip

(C1$'$) For each $x\in [2]^k_{H_1(i)}(2)$, there is a unique edge $e\in E_i'$ such that $e$ covers $x$.

\medskip

\noindent Now we divide the proof in three steps.

We first prove (c). By contradiction, suppose $I_{x}(e)\cap I_x(f)\ne\emptyset$. Pick  $i\in I_{x}(e)\cap I_x(f)$. Then $x\in [2]^k_{H_1(i)}(2)$, $\{e,f\}\subseteq L_i(H_2)$ and both $e$ and $f$ cover $x$, which implies that there is no edge in $E_i'$ covering $x$, contrary to (C1$'$).

The next thing is to prove (a). Suppose for the contrary that $|I_x(e)|\geq 2$. Pick distinct vertices $i,j\in I_x(e)$. Note that $E_i'\subseteq L_i(H_2)$ and $E_j'\subseteq L_j(H_2)$.  Combining (C1$'$) and (c), we get $e\in E_i'\cap E_j'$, contrary to (C2).

Finally, we prove (b). By contradiction, suppose that $I_x(e)\ne\emptyset$ and $I_y(e)\ne\emptyset$. Then $e=\{x,y\}$. Taking $i\in I_x(e)$ and $j\in I_y(e)$, we have $(y_{(i)},x_{(i)})=(1,2)$ and $(x_{(j)},y_{(j)})=(1,2)$, which implies that $i\neq j$. An argument similar to the one used in the previous step shows that  $e\in E_i'\cap  E_j'$, contrary to (C2).

The proof is complete.
$\qed$

As applications, we study   $K_{[k]}$-minimal graphs and $\overline K_{[k]}$-minimal graphs in the rest of this subsection. To inverstigate $K_{[k]}$-minimal graphs, we observe that
\begin{equation}\label{kk}
[2]^k_{K_{[k]}(i)}(2)=\{(2,\ldots,2)\},
\end{equation}
and then get the following example from Lemma~\ref{minimal bk2}.

\begin{eg}\label{uv}
 (i) Define $\mathbb U_k$ as the graph with the vertex set $[2]^k$ and the edge set
 $$
 E(\mathbb U_k)=\{(1,\ldots,1), (2,\ldots,2)\}.
 $$
 Then $\mathbb U_k$ is a $K_{[k]}$-minimal graph with size $1$.

 (ii) Define $\mathbb V_k$ as the graph with the vertex set $[2]^k$ and the edge set
 $$
 E(\mathbb V_k)=\bigcup_{i\in[k]}\{\{x,(2,\ldots,2)\}\mid x_{(i)}=1, x_{(t)}=2 \text{ for each }t\in[k]\setminus\{i\}\}.
 $$
 Then $\mathbb V_k$ is a $K_{[k]}$-minimal graph with size $k$.
\end{eg}

\begin{cor}
If $H_2$ is $K_{[k]}$-minimal, then $1\leq |E(H_2)|\leq k$,

  {\rm(i)} with the lower bound if and only if $H_2=\mathbb U_k$.

  {\rm(ii)} with the upper bound if and only if $H_2=\mathbb V_k$.
\end{cor}
\proof Noting that the degree of each vertex $i$ in $K_{[k]}$ is $k-1$, we have $1\leq |E(H_2)|\leq k$ by Theorem~\ref{enbounds12}.

(i) The ``if'' implication follows immediately from Example~\ref{uv}~(i), while the ``only if'' implication follows from (\ref{kk}) and Lemma~\ref{minimal bk2} (i).

(ii) We get the sufficiency from Example~\ref{uv}~(ii). To prove the necessity, suppose $|E(H_2)|=k$. Write $y=(2,\ldots,2)$. By (\ref{kk}), we have $J_y=[k]$ and $J_z=\emptyset$ for $z\in [2]^k\setminus\{y\}$.
 Noting that $|E(H_2)|$ attains the upper bound in Theorem~\ref{enbounds12}, we deduced from Corollary~\ref{tightcondition} (ii) and Lemma~\ref{minimal bk2} that $|I_y(e)|=1$ for $e\in E(H_2)$  and $\{I_y(e)\mid e\in E(H_2)\}$ is a partition of $[k]$.
For any $e\in E(H_2)$, write $e=\{x,y\}$. Choosing the unique vertex $i\in I_y(e)$, we have $x_{(i)}=1$ and $x_{(t)}=2$  for each $t\in[k]\setminus\{i\}$, as desired.
$\qed$

To investigate $\overline K_{[k]}$-minimal graphs, we observe that
\begin{equation}\label{nulli}
  [2]^k_{\overline K_{[k]}(i)}(2)=[2]^k_i(2),
\end{equation}
and  for each $x\in [2]^k$, we have
\begin{equation}\label{nullj}
  J_x=\{i\in[k]\mid x_{(i)}=2\}.
\end{equation}

\begin{eg}\label{rk}
  Define $\mathbb R_k$ as the graph with the vertex set $[2]^k$ and the edge set
 $$
 E(\mathbb R_k)=\{\{x,y\}\mid x,y\in[2]^k, x_{(i)}\neq y_{(i)} \text{ for each }i\in[k]\}.
 $$
 Then $\mathbb R_k$ is a $\overline K_{[k]}$-minimal graph with size $2^{k-1}$.
\end{eg}
\proof Noting that $E(\mathbb R_k)$ is a matching, we infer that $\mathbb R_k$ has size $2^{k-1}$ and each vertex is covered by a unique edge in $\mathbb R_k$.
Observing that $L_i(\mathbb R_k)=E(\mathbb R_k)$ for each $i\in[k]$,
we conclude that $\mathbb R_k$ is $\overline K_{[k]}$-minimal from Lemma~\ref{minimal bk} and (\ref{nulli}).
$\qed$

\begin{re}
  Actually, if $H_1$ has an isolated vertex,  then $\mathbb R_k$ is $H_1$-minimal. However, if $\overline K_{[k]}\prec H_1$, then $H_1\circ \mathbb R_k$ is not a minimal graph in $\mathcal B_k$.
\end{re}

Given a graph $G$, the {\em Cartesian product} of $s$ copies of $G$ is the graph $G^{\Box s}$ with the vertex set
$$
V(G^{\Box s})=\{(x_{(1)},\ldots,x_{(s)})\mid x_{(i)}\in V(G) \text{ for }i\in[s]\},
$$
where two vertices $x$ and $y$ are adjacent if and only if there exists an index $i\in[s]$ such that $\{x_{(i)},y_{(i)}\}\in E(G)$  and $x_{(j)}=y_{(j)}$ for all indices $j\in[s]\setminus\{i\}$.

\begin{eg}\label{p2k}
  Let $P_2$ be the graph with the vertex set $\{1,2\}$ and the edge set $\{\{1,2\}\}$.
 Then $P_2^{\Box k}$ is a $\overline K_{[k]}$-minimal graph with size $k\cdot 2^{k-1}$.
\end{eg}
\proof It is routine to verify that $|E(P_2^{\Box k})|=k\cdot 2^{k-1}$.
For an edge $\{x,y\}$ in $P_2^{\Box k}$ with $(x_{(i)},y_{(i)})=(1,2)$, by (\ref{nulli}), we have $I_x(\{x,y\})=\emptyset$ and $I_y(\{x,y\})=\{i\}$, and so $\tilde I_y(\{x,y\})=\{i\}$.
By Lemma~\ref{minimal bk2} and (\ref{nullj}), we conclude that $P_2^{\Box k}$ is $\overline K_{[k]}$-minimal.
$\qed$

\begin{cor}
If $H_2$ is $\overline K_{[k]}$-minimal, then $2^{k-1}\leq |E(H_2)|\leq k\cdot 2^{k-1}$,

  {\rm(i)} with the lower bound if and only if $H_2=\mathbb R_k$.

  {\rm(ii)} with the upper bound if and only if $H_2=P_2^{\Box k}$.
\end{cor}
\proof Since the degree of each vertex $i$ in $\overline K_{[k]}$ is $0$, by Theorem~\ref{enbounds12}, we have $2^{k-1}\leq |E(H_2)|\leq k\cdot 2^{k-1}$.

(i) The ``if'' implication follows immediately from Example~\ref{rk}, while the ``only if'' implication follows from (\ref{nulli}) and Corollary~\ref{tightcondition} (i).

(ii) We get the sufficiency from Example~\ref{p2k}. To prove the necessity,
suppose that the size of $H_2$ is $k\cdot 2^{k-1}$, as well as the upper bound in Theorem~\ref{enbounds12}  for $H_1=\overline K_{[k]}$. For $x\in[2]^k$,  write $\varepsilon(x)=\{e\in E(H_2)\mid I_x(e)\neq\emptyset\}$. It follows from (\ref{nullj}) and Lemma~\ref{minimal bk2} (i) that $\varepsilon(x)=\emptyset$ if and only if $x=(1,\ldots,1)$. Observing that $\tilde I_x(e)\subseteq I_x(e)$, we have derived from Lemma~\ref{minimal bk2} (ii) that
$$
E(H_2)=\bigcup_{x\in[2]^k}\{e\in E(H_2)\mid I_x(e)\neq\emptyset\}=\bigcup_{x\in[2]^k\setminus\{(1,\ldots,1)\}} \varepsilon(x).
$$
Pick any vertex $x\in[2]^k\setminus\{(1,\ldots,1)\}$. Combining Corollary~\ref{tightcondition}~(ii) and Lemma~\ref{minimal bk2}~(i), we get $|I_x(e)|=1$ for any $e\in\varepsilon(x)$ and $\{I_x(e)\mid e\in\varepsilon(x)\}$ is a partition of $J_x$.
For any edge $e=\{y,x\}\in\varepsilon(x)$, taking the unique vertex $i\in I_x(e)$, one has $(y_{(i)},x_{(i)})=(1,2)$ and $(y_{(t)},x_{(t)})=(2,2)$ for $t\in J_x\setminus\{i\}$, and further by Corollary~\ref{tightcondition}~(ii), we get $I_y(e)=\emptyset$, and so $(y_{(t)},x_{(t)})=(1,1)$ for $t\in [k]\setminus J_x$.
Hence, we have $\varepsilon(x)\in E(P_2^{\Box k})$, and so $E(H_2)\subseteq E(P_2^{\Box k})$. It follows from $|E(H_2)|=|E(P_2^{\Box k})|$ that $H_2=P_2^{\Box k}$, as desired.
$\qed$

\subsection{A poset on $\mathcal C_k$}

In this subsection, we focus on the poset $(\mathcal C_k,\preceq)$, and always suppose that $H_2$ is a graph with the vertex set $[3]^k$. We begin by constructing a graph on the set $[3]^k$.
Let $\Gamma_k$ be the graph with the vertex set $[3]^k$ and the edge set
$$
E(\Gamma_k)=\{\{x,y\}\mid x,y\in[3]^k,\; x\neq y,\; |x_{(i)}-y_{(i)}|\leq 1\text{ for each }i\in[k]\}.
$$

\begin{ob}\label{gammak}
  $\Gamma_k=span(\bigcup_{i\in [k]}(E(C^k_i)\cup E(D^k_i)))$.
\end{ob}

\begin{thm}\label{c adding edges}
  Suppose that $\overline K_{[k]}\circ H_2\in \mathcal C_k$.

  {\rm(i)} Then $H_2\preceq \Gamma_k$.

  {\rm(ii)} If $H_2\preceq H_2'\preceq \Gamma_k$, then $\overline K_{[k]}\circ H_2'\in\mathcal C_k$.
\end{thm}
\proof (i) Combining Proposition~\ref{lemmack} (i) and Observation~\ref{gammak}, we have  $E(H_2)\subseteq E(\Gamma_k)$, and so $H_2\preceq \Gamma_k$.

(ii) Observe that $E(H_2')\subseteq E(\Gamma_k)$, $M_i(H_2)\subseteq M_i(H_2')$ and $N_i(H_2)\subseteq N_i(H_2')$, where $M_i(H_2)$ and $N_i(H_2)$ are as refer to (\ref{mih2}) and (\ref{nih2}), respectively.
 For each $i\in[k]$, by Observation~\ref{gammak}, we have $M_i(\Gamma_k)=E(C^k_i)$ and $N_i(\Gamma_k)=E(D^k_i)$, which are  $[3]^k_i(2)$-covering and $S^k_i$-covering, respectively.
 Hence,  the desired result follows from Proposition~\ref{lemmack}.
$\qed$

\begin{cor}
For any graph $\overline K_{[k]}\circ H_2\in \mathcal C_k$, we have
$|E(H_2)|\leq \frac{7^k-3^k}{2},$
with equality if and only if $H_2=\Gamma_k$.
\end{cor}
\proof  Noting that $|x_{(i)}-y_{(i)}|\leq 1$ if and only if
$$
(x_{(i)},y_{(i)})\in\{(1,1),(1,2),(2,1),(2,2),(2,3),(3,2),(3,3)\},
$$
we have
$|\{(x,y)\mid x,y\in[3]^k,|x_{(i)}-y_{(i)}|\leq 1 \text{ for each }i\in[k]\}|= 7^k,$
which implies that $|E(\Gamma_k)|=\frac{7^k-3^k}{2}$. Hence, by Theorem~\ref{c adding edges}, we get the desired result.
$\qed$

\begin{cor}\label{unionc}
For $G_1,G_2\in\mathcal C_k$, we have $G_1\cup G_2\in \mathcal C_k$. In particular, we have
$$
\bigcup\limits_{G\in \mathcal C_k}G=\overline K_{[k]}\circ \Gamma_{k}\in\mathcal C_k.
$$
\end{cor}
\proof Write $G_1=\overline K_{[k]}\circ H_2$ and $G_2=\overline K_{[k]}\circ H_2'$. Then $G_1\cup G_2=\overline K_{[k]}\circ (H_2\cup H_2')$.
Hence, the desired result follows from Theorem~\ref{c adding edges}.
$\qed$

The following result is immediate from Corollary~\ref{unionc}.

\begin{prop}
 The poset $(\mathcal C_k,\preceq)$ is a finite join-semilattice with the maximum graph $\overline K_{[k]}\circ \Gamma_{k}$.
\end{prop}

A graph $H_2$ is {\em $k$-minimal} if $\overline K_{[k]}\circ H_2\in\mathcal C_k$ and $\overline K_{[k]}\circ H_2'\not\in\mathcal C_k$ for $H_2'\prec H_2$.

\begin{re}\label{remark5}
  A graph $H_2$ is $k$-minimal if and only if $\overline K_{[k]}\circ H_2$ is minimal in $\mathcal C_k$.
\end{re}

To study $k$-minimal graph, we give the following notation for convenience.

\begin{nota}
  Write $X=[3]^k_{[k]}(\{2,3\})$, $Y=[3]^k_{[k]}(\{1,3\})$ and
  $$
  \overline X=[3]^k\setminus X,\qquad \overline Y=[3]^k\setminus Y,\qquad Z=[3]^k\setminus (X\cup Y).
  $$
\end{nota}

We establish the bounds for the size of a $k$-minimal graph.

\begin{thm}\label{ceminimalbounds}
 If $H_2$ is $k$-minimal, then
  \begin{eqnarray*}
  \frac{3^k+1}{2}\leq |E(H_2)|\leq k\cdot(3^{k-1}+2^{k-1}).
  \end{eqnarray*}
\end{thm}
\proof For $i\in[k]$, denote by $M'_i$ (resp. $N'_i$) the set of edges $e$ in $M_i(H_2)$ (resp. $N_i(H_2)$) such that $M_i(H_2)\setminus\{e\}$ (resp. $N_i(H_2)\setminus\{e\}$) is not a $[3]^k_i(2)$-covering (resp. an $S^k_i$-covering), where $M_i(H_2)$ (resp. $N_i(H_2)$) is as refer to (\ref{mih2}) (resp. (\ref{nih2})). Then $|M_i'|\leq |[3]^k_i(2)|=3^{k-1}$ and $|N'_i|\leq |S^k_i|=2^{k-1}$. Since $H_2$ is $k$-minimal, it follows from Proposition~\ref{lemmack} that
$E(H_2)=\bigcup_{i\in[k]}(M_i'\cup N_i'),$
and so
\begin{equation}\label{m'n'}
  |E(H_2)|=|\bigcup_{i\in[k]}(M_i'\cup N_i')|\leq \sum_{i\in[k]}(|M_i'|+|N_i'|)\leq k\cdot(3^{k-1}+2^{k-1}).
\end{equation}

In the following, we prove the lower bound.
 Pick any vertex $x\in[3]^k$.
We shall find an edge $e_x$ in $H_2$ covering $x$ for the following two cases.

{\em Case 1.} $x=(3,\ldots,3)$. Then $x\in S^k_1$. Note that $\overline K_{[k]}\circ H_2\in\mathcal C_k$. By Proposition~\ref{lemmack}~(iii), there exists at least one edge $e_x$ in $H_2$ covering $x$.

{\em Case 2.} $x\in\overline Y$. Then there is an index $i\in[k]$ with $x_{(i)}=2$. By Proposition~\ref{lemmack}~(ii), there is an edge $e_x=\{x,y\}$ in $H_2$ such that $y_{(i)}=1$, and so $y\in \overline X$.

Observing that $X=(X\cap \overline Y)\cup\{(3,\ldots,3)\}$ and $Z\subseteq \overline Y$, we choose two families of edges from the above two cases:
$$
F_X=\{e_x\mid x\in X\}\qquad \text{and} \qquad F_Z=\{e_x\mid x\in Z\}.
$$
By Case 2, the edge $e_x$ for $x\in X\cap\overline Y$ does not cover any other vertices in $X$ except $x$, and so $|F_X|=|X|=2^k$. Noting that $Z\subseteq\overline X$, we have derived that all of the vertices covered by edges in $F_Z$ are in $\overline X$, which implies that $F_X\cap F_Z=\emptyset$. Since each edge covers two vertices, one gets
$2|F_Z|\geq |Z|=3^k-2^{k+1}+1$ by a short calculation.
Consequently, we have
\begin{equation}\label{xz}
|E(H_2)|\geq |F_X\cup F_Z|=|F_X|+|F_Z|\geq 2^k+\frac{3^k-2^{k+1}+1}{2}=\frac{3^k+1}{2},
\end{equation}
as desired.
$\qed$

\begin{re}
  If $\overline K_{[k]}\circ H_2$ is minimal in $\mathcal C_k$, then by the equation (\ref{eh1h2}), Remark~\ref{remark5} and Theorem~\ref{ceminimalbounds}, we have
  $$
   \frac{(2k+3)\cdot 3^{k-1}+1}{2}\leq |E(\overline K_{[k]}\circ H_2)|\leq 2k\cdot( 3^{k-1}+2^{k-2}).
  $$
\end{re}

 To characterize the lower bound in Theorem~\ref{ceminimalbounds}, we define a graph on $[3]^k$.

\begin{eg}\label{tk}
Define $\mathbb{T}_k$ as the graph with the vertex set $[3]^k$ and the edge set  $E(\mathbb T_k)=E_X\cup E_Z$, where
$$
  \begin{array}{rcl}
    E_{X}&=&\{\{x,y\}\mid x\in X,y\in [3]^k,\; x_{(i)}-y_{(i)}=1\text{ for each }i\in[k]\},\\
    E_{Z}&=&\{\{x,y\}\mid x,y\in Z,\; x_{(i)}=y_{(i)}=3\text{ or }\{x_{(i)},y_{(i)}\}=\{1,2\}\text{ for each }i\in[k]\}.
  \end{array}
  $$
  Then $\mathbb{T}_k$ is a $k$-minimal graph with size $\frac{3^k+1}{2}$.
\end{eg}
\proof We first compute $|E(\mathbb T_k)|$. Since there is a bijection from $X$ to $E_X$, one has $|E_X|=|X|=2^k$. Each vertex in $Z$ is covered by a unique edge in $E_Z$, so $E_Z$ is a matching in $Z$, which implies that
$$
|E_Z|= \frac{|Z|}{2}=\frac{3^k-2^{k+1}+1}{2}.
$$
Noting that $E_Z$ does not cover any vertex in $X$, we have $E_X\cap E_Z=\emptyset$. Hence, we get $|E(\mathbb T_k)|=|E_X|+|E_Z|=\frac{3^k+1}{2}$.

The next thing is to prove $\overline K_{[k]}\circ \mathbb T_k\in\mathcal C_k$. Take any index $i\in[k]$ and then pick any vertex $x\in[3]^k_i(2)\cup S_i^k$.
If $x\in X$, let $y$ be a vertex in $[3]^k$ such that $y_{(t)}=x_{(t)}-1$ for each $t\in[k]$, then $\{x,y\}\in E_X\cap E(C^k_i)$ or $E_X\cap E(D^k_i)$ according to $x\in [3]^k_i(2)$ or $S_i^k$.
If $x\in \overline X$, then $x\in Z$, and so $\{x,y\}\in E_Z\cap E(C^k_i)$, where  $y$ is a vertex in $Z$ such  that
$$
y_{(t)}=\left\{
\begin{array}{ll}
2,&\text{if }x_{(t)}=1,\\
1,&\text{if }x_{(t)}=2,\\
3,&\text{if }x_{(t)}=3.
\end{array}\right.
$$
Therefore, we have derived that $M_i(H_2)$, as refer to (\ref{mih2}), is a $[3]^k_i(2)$-covering and $N_i(H_2)$, as refer to (\ref{nih2}), is an $S^k_i$-covering. Observing that $E_X\cup E_Z\subseteq \bigcup_{i\in[k]}(E(C^k_i)\cup E(D^k_i))$, we get $\overline K_{[k]}\circ H_2\in\mathcal C_k$ from Proposition~\ref{lemmack}.

Finally, we show that $\mathbb T_k$ is $k$-minimal. By contradiction, if $\mathbb T_k$ is not $k$-minimal, then there is a graph $H_2\prec \mathbb T_k$ such that $\overline K_{[k]}\circ H_2\in\mathcal C_k$, and so $|E(H_2)|<|E(\mathbb T_k)|=\frac{3^k+1}{2}$, contrary to Theorem~\ref{ceminimalbounds}. This contradiction implies that $\mathbb T_k$ is minimal.

The proof is complete.
$\qed$

We use the following result to characterize the lower bound in Theorem~\ref{ceminimalbounds}.

\begin{cor}
  The lower bound in Theorem~{\rm \ref{ceminimalbounds}} is attained if and only if $H_2=\mathbb T_k$.
\end{cor}
\proof The ``if'' implication follows from Example~\ref{tk}. To obtain the ``only if'' implication,  let  $H_2$ be a $k$-minimal graph with size $\frac{3^k+1}{2}$.
With references to the proof of Theorem~\ref{ceminimalbounds}, we have derived from (\ref{xz}) that
$$
E(H_2)=F_X\cup F_Z, \quad |F_X|=2^k=|E_X|\quad \text{and}\quad |F_Z|=\frac{|Z|}{2}=|E_Z|.
$$
It follows from Proposition~\ref{lemmack} that $F_X=E_X$ and $F_Z=E_Z$, as desired.
$\qed$

To characterize the upper bound in Theorem~\ref{ceminimalbounds}, we give the following notation.

\begin{nota}
  For $i\in[k]$ and $x\in [3]_i^k(2)\cup S_i^k$, denote by $\epsilon_i(x)$ the set of edges $e=\{x,x'\}$ such that $x'_{(i)}=x_{(i)}-1$ and for each $t\in[k]\setminus\{i\}$, one of the following conditions holds.

  (i) If $x\in S_i^k$, then $x'_{(t)}=x_{(t)}$.

  (ii) If $x\in X\cap [3]_i^k(2)$, then
  $$
  x'_{(t)}=\left\{
  \begin{array}{ll}
  2\text{ or }3, &\text{if } x_{(t)}=2,\\
  3, &\text{if }x_{(t)}=3.
  \end{array}
  \right.
  $$

  (iii) If $x\in\overline X\cap [3]_i^k(2)$, then
  $$
  x'_{(t)}=\left\{
  \begin{array}{ll}
  1, &\text{if }x_{(t)}=1,\\
  2\text{ or }3, &\text{if } x_{(t)}=2\text{ or }3.
  \end{array}
  \right.
  $$
\end{nota}

\begin{lemma}\label{epsilonix}
  Let $i,j\in[k]$. Pick $x\in [3]_i^k(2)\cup S_i^k$ and $y\in [3]_j^k(2)\cup S_j^k$. If $(i,x)\neq (j,y)$, then $\epsilon_i(x)\cap \epsilon_j(y)=\emptyset$.
\end{lemma}
\proof Suppose for the contrary that $\epsilon_i(x)\cap \epsilon_j(y)\neq\emptyset$. Take $e\in\epsilon_i(x)\cap \epsilon_j(y)$.

{\em Case 1.} $i=j$. Then $x\neq y$ and $e=\{x,y\}$, which implies that $y_{(i)}=x_{(i)}-1$ and $x_{(i)}=y_{(i)}-1$, a contradiction.

{\em Case 2.} $i\neq j$ and $x=y$. Write $e=\{x,x'\}$. Then
$$
x'_{(i)}=x_{(i)}-1\quad \text{and}\quad x'_{(j)}=x_{(j)}-1.
$$
Note that $e\in\epsilon_i(x)$. If  $x\in S_i^k$ or $X\cap [3]_i^k(2)$, then $x'_j\geq x_j$, a contradiction. If $x\in \overline X\cap [3]_i^k(2)$, then $(x'_{(j)},x_{(j)})=(2,3)$, which implies that $x\in S^k_j$, and so $x\in X$, a contradiction.

{\em Case 3.} $i\neq j$ and $x\neq y$. Then $e=\{x,y\}$, and so
$$
y_{(i)}=x_{(i)}-1\quad \text{and}\quad y_{(j)}=x_{(j)}+1.
$$
If $x\in S_i^k$, then $y_{(j)}=x_{(j)}$, a contradiction.
Now suppose $x\in [3]_i^k(2)$. Then $y_{(i)}=1$ and $(x_{(j)},y_{(j)})=(2,3)$. Noting that $\{x,y\}\in\epsilon_j(y)$, we have
$y\in S^k_j$, and so $x_{(i)}=y_{(i)}$,  a contradiction.

We accomplish the proof.
$\qed$

\begin{eg}\label{qk}
  Define
  $$
  \mathcal Q_k=\{span(\bigcup_{i\in[k]}\;\bigcup_{x\in [3]_i^k(2)\cup S_i^k}\{e_i(x)\})\mid e_i(x)\in\epsilon_i(x)\}.
  $$
  Then each graph in $\mathcal Q_k$ is a $k$-minimal graph with size $k\cdot(3^{k-1}+2^{k-1})$.
\end{eg}
\proof Choose any graph $H_2\in\mathcal Q_k$. Observe that $\epsilon_i(x)\subseteq E(C^k_i)\cup E(D^k_i)$. Noting that $E(H_2)\cap\epsilon_i(x)$ has exactly one edge, as well as the unique edge covering $x$ in $M_i(H_2)$ or $N_i(H_2)$ according to $x\in [3]^k_i(2)$ or $S^k_i$, we deduce that $H_2$ is $k$-minimal from Proposition~\ref{lemmack} and Lemma~\ref{epsilonix}.
By Lemma~\ref{epsilonix} again, one has
$$
|E(H_2)|=k\cdot(|[3]^k_i(2)|+|S^k_i|)=k\cdot(3^{k-1}+2^{k-1}),
$$
as desired.
$\qed$

\begin{re}
  Let $P_3$ be the graph with the vertex set $\{1,2,3\}$ and the edge set $\{\{1,2\},\{2,3\}\}$.
Let $\mathbb Q_k$ be the graph obtained from the Cartesian product $P_3^{\Box k}$ by deleting the edges $\{x,x'\}$ such that there exist indices $i$ and $j$ in $[k]$ with $(x_i,x'_i)=(2,3)$ and $x_j=x'_j=1$. Then $\mathbb Q_k\in\mathcal Q_k$.
\end{re}

\medskip

We use the following result to characterize the upper bound in Theorem~\ref{ceminimalbounds}.

\begin{cor}
  The upper bound in Theorem~{\rm \ref{ceminimalbounds}} is attained if and only if $H_2\in\mathcal Q_k$.
\end{cor}
\proof We get the sufficiency from Example~\ref{qk}. To prove the necessity, let $H_2$ be a $k$-minimal graph with size $k\cdot(3^{k-1}+2^{k-1})$.
With references to the proof of Theorem~\ref{ceminimalbounds}, we have derived from (\ref{m'n'}) that
\begin{equation}\label{m'n'2}
 |M'_i|=3^{k-1}=|[3]^k_i(2)|,\qquad |N'_i|=2^{k-1}=|S^k_i|,
\end{equation}
and the following condition.

\medskip

 (A) All sets $M_i'$'s and $N_i'$'s are pairwise non-intersecting.

\medskip

Take any $x\in [3]^k_i(2)\cup S_i^k$. By (\ref{m'n'2}), there is a unique edge $f_i(x)$ covering $x$ in $M'_i$ or $N'_i$ according to $x\in [3]^k_i(2)$ or $S_i^k$. Note that $E(H_2)=\bigcup_{i\in[k]}(M_i'\cup N'_i)$.
To get the desired result, we only need to prove $f_i(x)\in\epsilon_i(x)$.

Suppose for the contrary that $f_i(x)\not\in\epsilon_i(x)$. It is routine to verify that there is an index $t\in[k]\setminus\{i\}$ such that either $(x,f_i(x))\in [3]^k_t(2)\times E(C^k_t)$ or $(x,f_i(x))\in S^k_t\times E(D^k_t)$.
Then $f_i(x)\in M_t(H_2)$ or $N_t(H_2)$ according to $x\in [3]^k_t(2)$ or $S^k_t$. Noting that $M_t'\subseteq M_t(H_2)$ and $N_t'\subseteq N_t(H_2)$, by the definitions of $M_t'$ and $N_t'$ in Theorem~\ref{ceminimalbounds}, we get $f_t(x)=f_i(x)$, which  contradicts (A).
This contradiction completes the proof.
$\qed$

\section{Perfectness-resolvable}

We begin this section by computing the diameters of graphs in $\mathcal B_k$ or $\mathcal C_k$.

\begin{prop}
  {\rm (i)} The diameter of any graph  in $\mathcal B_k$ is $2$ or $3$.

  {\rm (ii)} The diameter of any graph  in $\mathcal C_k$ is $3$, $4$ or $5$.
\end{prop}
\proof (i) Choose any graph $H\in\mathcal B_k$. Pick $i,j\in[k]$. By Lemma~\ref{verifyb} (i), we have $d(i,z)=z_{(i)}$ for each $z\in[2]^k$, and so
$$
d(i,j)\leq d(i,(1,\ldots,1))+d(j,(1,\ldots,1))= 1+1=2.
$$
For distinct vertices $x,y\in[2]^k$, there exists a vertex $t\in[k]$ such that $x_{(t)}=1$ or $y_{(t)}=1$. Hence, we have
$$
d(x,y)\leq d(t,x)+d(t,y)\leq 1+2=3.
$$
Since $d(i,(2,\ldots,2))=2$, the diameter of $H$ is $2$ or $3$.

(ii) An argument similar to the proof of (i) shows that (ii) holds.
$\qed$

\begin{re}
  (i) The diameters of $K_{[k]}\circ K_{[2]^k}$ and $K_{[k]}\circ \mathbb U_2$ are $2$ and $3$, respectively.

  (ii) The diameters of $\overline K_{[k]}\circ \Gamma_k$, $\overline K_{[k]}\circ \mathbb Q_2$ and $\overline K_{[k]}\circ \mathbb T_2$ are $3$, $4$ and $5$, respectively.
\end{re}

Given a  graph $G$, the {\em metric dimension} of $G$, denoted by $\dim(G)$, is the minimum cardinality of a resolving set of $G$. A {\em metric basis} of $G$ is a resolving set of $G$ with cardinality $\dim(G)$. A metric basis of $G$ is {\em perfect} if it is a completeness-resolving set. We say that $G$ is {\em perfectness-resolvable} if it admits a perfect metric basis. Clearly, a perfectness-resolvable graph is completeness-resolvable.

\begin{ob}
  (i) All paths are perfectness-resolvable.

  (ii) A graph in $\mathcal K$ is perfectness-resolvable if and only if it is complete.
\end{ob}

\begin{prop}
  {\rm(i)} Let $G$ be a graph in $\mathcal B_k$. If the diameter of $G$ is $2$, then $G$ is perfectness-resolvable.

  {\rm(ii)} Let $G$ be a graph in $\mathcal C_k$. If the diameter of $G$ is $3$, then $G$ is perfectness-resolvable.
\end{prop}
\proof By \cite[Theorem 1]{chartrandejo}, we have $|V(G)|\leq\dim(G)+d^{\dim(G)}$, where $d$ is the diameter of $G$. Note that $|V(G)|$ is equal to $k+2^k$ or $k+3^k$ according to $G\in\mathcal B_k$ or $\mathcal C_k$. Hence, one has $k\leq\dim(G)$ if the condition in (i) or (ii) holds. Since $[k]$ is a completeness-resolving set of graphs in $\mathcal B_k\cup\mathcal C_k$, the two desired results follow.
$\qed$

We conclude the paper by raising the following problem.

\medskip

\noindent {\bf Problem~2.} Which graphs are perfectness-resolvable?

\section*{Acknowledgements}

Feng was supported by the National Natural Science Foundation of China (11701281), the Natural Science Foundation of Jiangsu Province (BK20170817) and the Grant of China Postdoctoral Science Foundation. Ma was supported
by the National Natural Science Foundation of China (11801441) and the Young Talent fund of University Association for Science and Technology in Shanxi, China (20190507). Xu was supported by the National Natural Science Foundation of China (61673218) and the Postgraduate Research \& Practice
Innovation Program of Jiangsu Province (KYCX18\_0376
and KYCX19\_0253).

\end{CJK*}

\end{document}